# CONSISTENCY OF A RECURSIVE ESTIMATE OF MIXING DISTRIBUTIONS

By Surya T. Tokdar, Ryan Martin and Jayanta K. Ghosh

*Carnegie Mellon University, Purdue University and Purdue University, Indian Statistical Institute*

Mixture models have received considerable attention recently and Newton [*Sankhyā Ser. A* **64** (2002) 306–322] proposed a fast recursive algorithm for estimating a mixing distribution. We prove almost sure consistency of this recursive estimate in the weak topology under mild conditions on the family of densities being mixed. This recursive estimate depends on the data ordering and a permutation-invariant modification is proposed, which is an average of the original over permutations of the data sequence. A Rao–Blackwell argument is used to prove consistency in probability of this alternative estimate. Several simulations are presented, comparing the finite-sample performance of the recursive estimate and a Monte Carlo approximation to the permutation-invariant alternative along with that of the nonparametric maximum likelihood estimate and a nonparametric Bayes estimate.

**1. Introduction.** Mixture distributions have played a key role in modeling data that reflect population heterogeneity, contain indirect observations or involve latent variables. In recent years, these models have been widely used in genetics, bioinformatics, proteomics, computer vision, speech analysis and a host of other research areas; see, for example, [1, 5, 7, 21, 25, 27, 31]. Fitting a mixture model has been made easy by the advent of computational techniques such as the Expectation Maximization (EM) and the Markov Chain Monte Carlo (MCMC) algorithms. Recovering the underlying mixing distribution, however, continues to pose a serious challenge.

Newton, et al. [22, 23, 24] introduced a fast, recursive algorithm for estimating a mixing density when a finite sample is available from the corresponding mixture model. Suppose $X_1, \ldots, X_n$ are independently distributed

(iid) according to the density

$$p(x) = \int_{\Theta} p(x|\theta) F(d\theta), \tag{1}$$

where $p(x|\theta)$ is a known sampling density, with respect to a dominating $\sigma$-finite measure $\nu$ on $\mathcal{X}$, parametrized by $\theta \in \Theta$. Assume also that the mixing distribution $F$ is absolutely continuous with respect to some $\sigma$-finite measure $\mu$ on $\Theta$. Newton [22] proposed to estimate $f = dF/d\mu$ as follows.

RECURSIVE ALGORITHM. Fix an initial estimate $f_0$ and a sequence of weights $w_1, w_2, \ldots \in (0, 1)$. Given i.i.d. observations $X_1, \ldots, X_n$ from the mixture density $p(x)$ in (1), compute

$$f_i(\theta) = (1 - w_i) f_{i-1}(\theta) + w_i \frac{p(X_i|\theta) f_{i-1}(\theta)}{\int_{\Theta} p(X_i|\theta') f_{i-1}(\theta') \mu(d\theta')}, \qquad \theta \in \Theta \tag{2}$$

for $i = 1, \ldots, n$ and produce $f_n$ as the final estimate.

This method of estimating $f$ has a number of advantages over the existing mainstream methods found in the literature. To begin with, it is rather straightforward to accommodate prior information regarding support and continuity properties of $f$ through those of $f_0$. For example, if $p(x_i|\theta) > 0$ for all $\theta \in \Theta$, then $\mathrm{supp}(f_i) = \mathrm{supp}(f_{i-1})$. Therefore, by choosing $f_0$ appropriately one can ensure that $f_n$ has the same support as the target $f$. This flexibility is not offered by the method of nonparametric maximum likelihood (NPML) estimation of $F$; cf. Laird [16] and Lindsay [18] which produces an estimate supported on at most $n$ points. It is also evident that the recursive algorithm above applies to any arbitrary sampling density $p(x|\theta)$, making this method more general than deconvolution methods which deal exclusively with sampling densities of the type $p(x|\theta) = \varphi(x - \theta)$ for some density $\varphi$. However, a lot is known about deconvolution methods; see, for example [9, 30, 35].

The flexibility associated with $f_n$ resembles those found in a Bayesian framework. Indeed, for $n = 1$, the estimate $f_n$ is precisely the posterior mean of $f$ under the Bayesian formulation that a priori $f$ follows a Dirichlet process (DP) distribution [8, 10] with base measure $f_0$ and precision $1/w_1 - 1$. Newton's original motivation for the recursive algorithm was based on this fact [23], though this analogy breaks down for $n > 1$. In particular, $f_n$, for $n > 1$, depends on the particular order in which $X_i$'s enter the recursion and hence is not a function of the sufficient statistic $(X_{(1)}, \ldots, X_{(n)})$—the order statistics. Consequently, $f_n$ cannot equal any posterior quantity. To further distinguish the two estimates, computation of a nonparametric Bayes estimate $f_{\mathrm{NPB}}$ based on the DP prior requires a rather elaborate Monte Carlo procedure, while the recursive estimate can be computed many times faster.

It was hoped that $f_n$ would serve as a computationally efficient approximation to $f_{\mathrm{NPB}}$.

It is rather difficult to study the asymptotic properties of $f_n$ since it is not a Bayesian quantity, does not seem to optimize a criterion function such as the log likelihood and cannot be written as a linear estimator [13]. Nevertheless, empirical studies carried out by Newton [22] and Ghosh and Tokdar [12] clearly demonstrated good performance of this estimate. Newton [22] also presented a proof of convergence of $f_n$ as $n \to \infty$ based on the theory of inhomogeneous Markov chains. Unfortunately, this proof had a gap [12]. Ghosh and Tokdar [12] used a novel martingale based argument to show consistency under the same conditions as in Newton [22]. A slightly stronger result has recently been derived by Martin and Ghosh [20] using a stochastic approximation representation of the algorithm. The conditions required in these papers are somewhat restrictive, particularly in requiring $\Theta$ to be a known finite set. Ghosh and Tokdar [12] further require $p(x|\theta)$ to be bounded away from zero on $\mathcal{X} \times \Theta$, while Martin and Ghosh [20] make the weaker assumption that $p(\cdot|\theta) > 0$ $\nu$-a.e. for each $\theta$.

In this paper we show consistency of $f_n$ in the weak topology under quite general conditions. Our assumptions are slightly stronger than those needed in Kiefer and Wolfowitz [14] or Leroux [17] to prove consistency of the NPML estimate of $F$. In particular, $\Theta$ is not required to be finite and the sampling density is allowed to decay to zero. We use a major extension of the basic martingale argument in [12] to show that the marginal density

$$(3) \qquad p_n(x) = \int p(x|\theta) f_n(\theta) \mu(d\theta)$$

almost surely converges to $p(x)$ in the $L_1$ topology. This then leads to almost sure weak convergence of $f_n$ to $f$. This latter result applies to any arbitrary method of estimation—for example, it can be used to show weak consistency of the posterior mean of the mixing distribution under the DP formulation. This result holds even when $\Theta$ is noncompact, as long as $\theta \mapsto p(x|\theta)$ vanishes at the boundary in a certain near-uniform sense. The main martingale argument, too, does not explicitly require much structural assumption on $\Theta$, but assumption A5 (see Section 2) would be difficult to verify without compactness.

Despite this asymptotic justification, the dependence of $f_n$ on the order of the observations could be a cause of concern in application, especially when $n$ is not very large. In some cases a particular ordering can be justified by problem specific considerations. For example, "sparseness" assumptions (i.e., that only a small percentage of the observations come from the nonnull component of the mixture) led Bogdan, et al. [3] to arrange the observations in the ascending order of their magnitude while estimating a mixing density underlying a multiple testing problem. In the absence of such justification a permutation invariant version of $f_n$ may be desirable. Newton [22] recommends calculating the average over a large number of random permutations which can be seen as a Monte Carlo approximation to

$$\bar{f}_n = \frac{1}{n!} \sum_{s \in S_n} f_{n,s}, \tag{4}$$

where $S_n$ is the permutation group on $\{1, \ldots, n\}$ and $f_{n,s}$, for $s \in S_n$, represents the estimate $f_n$ with the observations arranged as $X_{s(1)}, \ldots, X_{s(n)}$.

In Section 3 we show that $\bar{f}_n$ provides a Rao–Blackwellization of $f_n$ and satisfies $\mathbb{E}d(f, \bar{f}_n) \leq \mathbb{E}d(f, f_n)$ for many standard divergence measures $d$. This property is then exploited to show that in the weak topology, $\bar{f}_n \to f$ in probability. Section 4 presents a simulation study of finite sample performance of $f_n$ and $\hat{f}_n$—a Monte Carlo approximation to $\bar{f}_n$. It is demonstrated that $\hat{f}_n$, which requires more computing time than $f_n$, is still faster and more accurate than other existing methods, such as the NPML estimate or a NP Bayes estimate. Finally, in Section 5 we give some concluding remarks.

**2. Consistency of Newton's estimate.** For the remainder of the paper we consider the following assumptions:

A1. $\sum_{i=1}^{\infty} w_i = \infty$ and $\sum_{i=1}^{\infty} w_i^2 < \infty$.
A2. The map $F \mapsto \int p(x|\theta)F(d\theta)$ is injective; that is, the mixing distribution $F$ is identifiable from the mixture $\int p(x|\theta)F(d\theta)$.
A3. For each $x \in \mathcal{X}$, the map $\theta \mapsto p(x|\theta)$ is bounded and continuous.
A4. For any $\varepsilon > 0$ and any compact $\mathcal{X}_0 \subset \mathcal{X}$, there exists a compact $\Theta_0 \subset \Theta$ such that $\int_{\mathcal{X}_0} p(x|\theta)\nu(dx) < \varepsilon$ for all $\theta \notin \Theta_0$.
A5. There exists a constant $B < \infty$ such that, for every $\theta_1, \theta_2, \theta_3 \in \Theta$

$$\int_{\mathcal{X}} \left( \frac{p(x|\theta_1)}{p(x|\theta_2)} \right)^2 p(x|\theta_3)\nu(dx) < B.$$

The first condition on the weights $w_i$ is necessary for $f_n$ to outgrow the influence of the initial guess $f_0$. At the same time, the weights need to decay to zero to allow for accumulation of information. The square summability condition ensures a certain rate for this decay, suitable for the Taylor approximation approach taken here. The identifiability condition A2, necessary for any estimation of mixture densities, is shown in Teicher [33] to be satisfied by many sampling densities of interest; for example:

- Normal with mean $\theta$ and fixed variance $\sigma^2 > 0$,
- Gamma with rate $\theta$ and fixed shape $\alpha > 0$,
- Poisson with mean $\theta$.

Each of these densities satisfy the boundedness conditions A3 as well as the decay property A4. These also satisfy the square integrability condition A5 when $\Theta$ is a compact interval.

Let $K_n = \int f \log(f/f_n) \, d\mu$ and $K_n^* = \int p \log(p/p_n) \, d\nu$ denote the error measures according to the Kullback–Leibler (KL) divergence, where $f_n$ and $p_n$ are defined in (2) and (3), respectively. On application of a telescoping sum, it follows easily from the recursive definition of the estimates $f_i$ that

$$K_n - K_0 = -\sum_{i=1}^{n} \int_{\Theta} \log\left[1 + w_i\left(\frac{p(X_i|\theta)}{p_{i-1}(X_i)} - 1\right)\right] f(\theta)\mu(d\theta).$$

Write $\log(1+x) = x - x^2 R(x)$, $x > -1$, where the remainder term satisfies $0 \le R(x) \le \max(1, (1+x)^{-2})/2$. Then

$$
\begin{aligned}
K_n - K_0 &= \sum_{i=1}^{n}\left[w_i\left(1 - \frac{p(X_i)}{p_{i-1}(X_i)}\right) + \int_{\Theta} R_i(X_i,\theta) f(\theta)\mu(d\theta)\right] \\
&= \sum_{i=1}^{n} w_i V_i - \sum_{i=1}^{n} w_i M_i + \sum_{i=1}^{n} E_i,
\end{aligned}
$$

(5)

where

$$R_i(x,\theta) = w_i^2\left(\frac{p(x|\theta)}{p_{i-1}(x)} - 1\right)^2 R\left(w_i\left(\frac{p(x|\theta)}{p_{i-1}(x)} - 1\right)\right),$$

$$M_i = -\mathbb{E}\left[1 - \frac{p(X_i)}{p_{i-1}(X_i)}\Big|\mathscr{F}_{i-1}\right]$$

$$= \int_{\mathcal{X}}\left(\frac{p(x)}{p_{i-1}(x)} - 1\right)p(x)\nu(dx),$$

$$V_i = 1 - \frac{p(X_i)}{p_{i-1}(X_i)} + M_i,$$

$$E_i = \int_{\Theta} R_i(X_i,\theta) f(\theta)\mu(d\theta).$$

In the following we prove that each of $\sum_{i=1}^{\infty} w_i V_i$, $\sum_{i=1}^{\infty} E_i$ and $\sum_{i=1}^{\infty} w_i M_i$ is finite with probability 1.

Let $\mathscr{F}_i = \sigma(X_1, \ldots, X_i)$ be the $\sigma$-algebra generated by the first $i$ observations. By definition of $V_i$, $S_n = \sum_{i=1}^{n} w_i V_i$ forms a zero mean martingale sequence with respect to $\mathscr{F}_n$. Moreover by assumption A5,

$$
\begin{aligned}
\mathbb{E}(S_n^2) &= \mathbb{E}\left\{\sum_{i=1}^{n} w_i^2 \mathbb{E}[V_i^2|\mathscr{F}_{i-1}]\right\} \\
&\le \mathbb{E}\left\{\sum_{i=1}^{n} w_i^2 \mathbb{E}\left[\left(1 - \frac{p(X_i)}{p_{i-1}(X_i)}\right)^2\Big|\mathscr{F}_{i-1}\right]\right\} \\
&\le 2(1+B)\sum_{i=1}^{n} w_i^2.
\end{aligned}
$$

Since $\sum_{i=1}^{\infty} w_i^2 < \infty$ by A1, $\mathbb{E}(S_n^2)$ is uniformly bounded in $n$. Therefore, by the martingale convergence theorem [6], $S_n$ almost surely converges to a random variable $S_\infty$ with $\mathbb{E}(S_\infty) < \infty$.

Let $T_n = \sum_{i=1}^n E_i$ and let $T_\infty = \lim T_n$, which always exists (since $E_i$'s are nonnegative) but may equal infinity. Notice that for $u > 0$ and $v \in (0,1)$,

$$(u-1)^2 \max\{1, (1 + v(u-1))^{-2}\} \leq \max\{(u-1)^2, (1/u-1)^2\}.$$

Therefore,

$$R_i(x, \theta) \leq \frac{w_i^2}{2} \left( \frac{p(x|\theta)}{p_{i-1}(x)} - 1 \right)^2 \max\left\{1, \left(1 + w_i \left( \frac{p(x|\theta)}{p_{i-1}(x)} - 1 \right) \right)^{-2} \right\}$$

$$\leq \frac{w_i^2}{2} \max\left\{ \left(1 - \frac{p(x|\theta)}{p_{i-1}(x)} \right)^2, \left(1 - \frac{p_{i-1}(x)}{p(x|\theta)} \right)^2 \right\}$$

and hence, by assumption A5,

$$\mathbb{E}[E_i | \mathscr{F}_{i-1}] = \int_\Theta \mathbb{E}[R_i(X_i, \theta) | \mathscr{F}_{i-1}] f(\theta) \mu(d\theta) \leq w_i^2 (1 + B).$$

By Fatou's lemma and assumption A1,

$$\mathbb{E}(T_\infty) \leq \liminf_n \mathbb{E}(T_n) = \liminf_n \mathbb{E}\left( \sum_{i=1}^n \mathbb{E}[E_i | \mathscr{F}_{i-1}] \right)$$

$$\leq \liminf_n (1 + B) \sum_{i=1}^n w_i^2 < \infty,$$

which proves $T_\infty$ is finite with probability 1.

Now rearrange the terms in (5) and use nonnegativity of $K_n$ to get

$$\sum_{i=1}^n w_i M_i \leq K_0 + S_n + T_n. \tag{6}$$

It follows from the inequality $\log y \leq y - 1$ that $M_i \geq K_i^* \geq 0$. Therefore $\sum_{i=1}^\infty w_i M_i$ exists but could be infinite. However, equation (6) implies

$$\sum_{i=1}^\infty w_i M_i \leq K_0 + S_\infty + T_\infty < \infty \qquad \text{a.s.} \tag{7}$$

The almost sure finiteness of the three series in (5) leads to the following important result.

THEOREM 1. *Under* A1 *and* A5, $K_n \to K_\infty$ *a.s. for some random variable* $K_\infty$. *Moreover,* $K_n^* \to 0$ *a.s. on a (random) subsequence.*

PROOF. The first assertion is a simple consequence of the finiteness of the three series. The second observation follows since $\sum_{i=1}^{\infty} w_i K_i^* < \infty$ almost surely while $\sum_{i=1}^{\infty} w_i = \infty$. □

Next, define the quantities

$$g_{i,x}(\theta) = \frac{p(x|\theta)f_{i-1}(\theta)}{p_{i-1}(x)} \quad \text{and} \quad h_{i,x'}(x) = \int_{\Theta} p(x|\theta)g_{i,x'}(\theta)\mu(d\theta),$$

so that the recursive updates $f_{i-1} \mapsto f_i$ and $p_{i-1} \mapsto p_i$ are, respectively,

$$f_i(\theta) = (1 - w_i)f_{i-1}(\theta) + w_i g_{i,X_i}(\theta),$$
$$p_i(x) = (1 - w_i)p_{i-1}(x) + w_i h_{i,X_i}(x).$$

Therefore, as in the case of $K_n$, we could write

$$K_n^* - K_0^* = -\sum_{i=1}^{n} \int_{\mathcal{X}} p(x) \log\left[1 + w_i\left(\frac{h_{i,X_i}(x)}{p_{i-1}(x)} - 1\right)\right]\nu(dx)$$

$$= \sum_{i=1}^{n} \int_{\mathcal{X}} p(x)\left[w_i\left(1 - \frac{h_{i,X_i}(x)}{p_{i-1}(x)}\right) + R_i^*(X_i, x)\right]\nu(dx)$$

$$= \sum_{i=1}^{n} w_i V_i^* - \sum_{i=1}^{n} w_i M_i^* + \sum_{i=1}^{n} E_i^*,$$

where

$$R_i^*(x', x) = w_i^2\left(\frac{h_{i,x'}(x)}{p_{i-1}(x)} - 1\right)^2 R\left(w_i\left[\frac{h_{i,x'}(x)}{p_{i-1}(x)} - 1\right]\right),$$

$$E_i^* = \int_{\mathcal{X}} R_i^*(X_i, x)p(x)\nu(dx),$$

$$M_i^* = -\mathbb{E}\left[1 - \int_{\mathcal{X}} \frac{h_{i,X_i}(x)}{p_{i-1}(x)}p(x)\nu(dx)\,\bigg|\,\mathscr{F}_{i-1}\right]$$

$$= \int_{\mathcal{X}}\int_{\mathcal{X}} \frac{h_{i,x'}(x)}{p_{i-1}(x)}p(x)p(x')\nu(dx)\nu(dx') - 1,$$

$$V_i^* = 1 - \int_{\mathcal{X}} \frac{h_{i,X_i}(x)}{p_{i-1}(x)}p(x)\nu(dx) + M_i^*.$$

Proceeding as in the lead up to Theorem 1, it can be shown that each of $\sum_{i=1}^{\infty} w_i V_i^*$, $\sum_{i=1}^{\infty} w_i M_i^*$ and $\sum_{i=1}^{\infty} E_i^*$ is finite almost surely. The required nonnegativity of $M_i^*$ is established using Jensen's inequality as follows:

$$M_i^* = \int_{\mathcal{X}}\int_{\mathcal{X}} \frac{\int_{\Theta} p(x|\theta)p(x'|\theta)f_{i-1}(\theta)\mu(d\theta)}{p_{i-1}(x)p_{i-1}(x')}p(x)p(x')\nu(dx)\nu(dx') - 1$$

$$= \int_{\Theta} \left\{ \int_{\mathcal{X}} \frac{p(x|\theta)}{p_{i-1}(x)} p(x) \nu(dx) \right\}^2 f_{i-1}(\theta) \mu(d\theta) - 1$$

$$\geq \left\{ \int_{\Theta} \int_{\mathcal{X}} \frac{p(x|\theta)}{p_{i-1}(x)} p(x) \nu(dx) f_{i-1}(\theta) \mu(d\theta) \right\}^2 - 1$$

$$= 0.$$

From this we conclude the following.

THEOREM 2. *Under* A1 *and* A5, $K_n^* \to 0$ *a.s.*

PROOF. It follows from the above discussion that $K_n^* \to K_\infty^*$, almost surely, for some random variable $K_\infty^*$. But we know from Theorem 1 that $K_n^* \to 0$ a.s. on a subsequence. These two together imply $K_\infty^* = 0$ a.s. □

It follows from the above theorem that $p_n \to p$ in the $L_1$ topology (also in the topology of the Hellinger metric). We next show that under identifiability, $L_1$ consistency of $p_n$ implies weak consistency of $f_n$ via a tightness argument. Since this result requires no assumption on the construction of $f_n$, it is presented in the next theorem in more generality than required at the moment. We will see some other use of it in the sequel.

THEOREM 3. *Let* $\widetilde{F}$ *and* $\widetilde{F}_n$ *be probability measures on* $\Theta$ *with respective mixture densities* $\tilde{p}(x) = \int p(x|\theta)\widetilde{F}(d\theta)$ *and* $\tilde{p}_n(x) = \int p(x|\theta)\widetilde{F}_n(d\theta)$. *Suppose* $\tilde{p}_n \to \tilde{p}$ *in* $L_1$. *Then, under* A2–A4, $\widetilde{F}_n \to \widetilde{F}$ *in the weak topology.*

PROOF. We first show that $\widetilde{F}_n$ forms a tight sequence. Fix any $\varepsilon > 0$. It suffices to show existence of a compact $\Theta_0 \subset \Theta$ such that $\widetilde{F}_n(\Theta_0) > 1 - \varepsilon$ for sufficiently large $n$. Take any compact $\mathcal{X}_0 \subset \mathcal{X}$ such that $\int_{\mathcal{X}_0} \tilde{p} \, d\nu > 1 - \varepsilon/2$. By A4, there exists a compact $\Theta_0$ such that $\int_{\mathcal{X}_0} p(x|\theta)\nu(dx) < \varepsilon/2$ for all $\theta \notin \Theta_0$. Now apply the $L_1$ convergence of $\tilde{p}_n$ to $\tilde{p}$ to conclude

$$1 - \frac{\varepsilon}{2} < \int_{\mathcal{X}_0} \tilde{p} \, d\nu = \lim_{n \to \infty} \int_{\mathcal{X}_0} \tilde{p}_n \, d\nu \leq \liminf_{n \to \infty} \left\{ \widetilde{F}_n(\Theta_0) + \frac{\varepsilon}{2} \widetilde{F}_n(\Theta_0^c) \right\}.$$

Thus, $\widetilde{F}_n$ is tight and the final assertion will follow once we show every weakly convergent subsequence $\widetilde{F}_{n(k)}$ converges to $\widetilde{F}$. Now, if $\widetilde{F}_{n(k)} \to \widetilde{F}^*$ for some $\widetilde{F}^* \in \mathscr{P}(\Theta)$ then, by assumption A3, $\tilde{p}_{n(k)} \to \tilde{p}^*$ pointwise and hence in the $L_1$ topology (via Scheffé's theorem), where $\tilde{p}^*(x) = \int p(x|\theta)\widetilde{F}^*(d\theta)$. Therefore $\tilde{p}^* = \tilde{p}$, which, under A2, implies $\widetilde{F}^* = \widetilde{F}$. □

The following result precisely states what we have already proved regarding consistency of $f_n$.

COROLLARY 4. *Under* A1–A5, *the estimate* $f_n$ *obtained from* (2) *converges almost surely to* $f$ *in the weak topology.*

**3. Averaging over permutations.** It is easy to see that the permutation averaged estimate $\bar{f}_n$ can be written as $\bar{f}_n = \mathbb{E}[f_n \mid X_{(1)}, \ldots, X_{(n)}]$. Let $\bar{p}_n$ denote the corresponding mixture density

$$\bar{p}_n(x) = \int_\Theta p(x|\theta)\bar{f}_n(\theta)\mu(d\theta).$$

Then $\bar{p}_n$ also satisfies $\bar{p}_n = \mathbb{E}[p_n|X_{(1)}, \ldots, X_{(n)}]$. Therefore $\bar{f}_n$ and $\bar{p}_n$ produce a Rao–Blackwellization of $f_n$ and $p_n$, respectively, by making these functions of the sufficient statistic—the order statistics. As one might guess, this results in a smaller expected error in estimation, when error is measured by a divergence $d$ that is convex in the estimate

$$\begin{aligned}
\mathbb{E}\,d(f, \bar{f}_n) &= \mathbb{E}d(f, \mathbb{E}[f_n|X_{(1)}, \ldots, X_{(n)}]) \\
&\leq \mathbb{E}\{\mathbb{E}[d(f, f_n)|X_{(1)}, \ldots, X_{(n)}]\} \\
&= \mathbb{E}\,d(f, f_n)
\end{aligned}$$

and similarly, $\mathbb{E}d(p, \bar{p}_n) \leq \mathbb{E}d(p, p_n)$. Examples of such divergence measures $d$ include the KL divergence and the $L_1$ distance.

We next show that the above result leads to weak convergence of $\bar{f}_n$ to $f$. However, we prove convergence only in probability and not almost surely. Recall that $Y_n \to Y$ in probability if and only if every subsequence $n_k$ contains a further subsequence $n_{k(l)}$ such that $Y_{n_{k(l)}} \to Y$ a.s., whenever the underlying topology is metrizable.

THEOREM 5. *Under* A1–A5, *$\bar{f}_n$ converges weakly to $f$ in probability.*

PROOF. From Theorem 2 it follows that $\|p - p_n\|_1 \to 0$ a.s. Since the $L_1$ distance is bounded by 2, it follows by the dominated convergence theorem that $\mathbb{E}\|p - p_n\|_1 \to 0$. Rao–Blackwellization implies $\mathbb{E}\|p - \bar{p}_n\|_1 \to 0$ and, hence, $\|p - \bar{p}_n\|_1 \to 0$ in probability. Take an arbitrary subsequence $n_k$. It must contain a further subsequence $n_{k(l)}$ such that $\|p - \bar{p}_{n_{k(l)}}\|_1 \to 0$ a.s. Then Theorem 3 implies that $\bar{f}_{n_{k(l)}} \to f$ a.s. in the weak topology. The assertion follows since the weak topology is metrizable. □

REMARK 6. Even for moderate $n$, there are too many permutations to compute $\bar{p}_n$ exactly, so the Monte Carlo estimate $\hat{p}_n$ is used as a numerical approximation. Therefore, what we can conclude from Theorem 5 is a sort of *practical* consistency of $\hat{p}_n$; that is, for large $n$ and sufficiently many random permutations, $\hat{p}_n \approx \bar{p}_n$ and $\bar{p}_n \approx p$, which implies $\hat{p}_n \approx p$.

**4. Simulations.** The numerical results in [12, 22] show that $f_n$ performs well in a variety of problems. In the following subsections we compare, more extensively, the performance of the recursive estimate (RE) and the recursive estimate averaged over permutations (PARE), starting with initial guess $f_0$, with that of several popular competitors, namely, the nonparametric maximum likelihood estimate (MLE) and the nonparametric Bayes estimate (NPB) based on a Dirichlet process prior with base measure $f_0$ and precision constant set to 1. While RE and PARE are easy to compute, computation of MLE and NPB is nontrivial. For the MLE, we implement an efficient new algorithm of Wang [34]. To find NPB, we employ a new importance sampling method, based on a collapsing of the Pólya Urn scheme; see the Appendix. We set the following simulation parameters:

- $T = 100$ samples of size $n = 200$ are taken from the model.
- For PARE, 100 random permutations of the data are selected.
- For RE and PARE, the weights satisfy $w_i = (i+1)^{-1}$.
- For NPB, $R = 10{,}000$ importance samples are used; see (12).

The efficiency of our NPB algorithm is measured by the *effective sample size* (ESS) [19]. For an importance sample of size $R$ the ESS, given by

$$\text{ESS} = \frac{R}{1 + \text{var}\{\omega_1^*, \dots, \omega_R^*\}},$$

estimates the size of an "equivalent" i.i.d. sample from the posterior distribution of $f$, where $\omega_r^*$ is a *normalization* of the weight $\omega_r$ in (11).

REMARK 7. Estimation of a mixing distribution in the Dirichlet process mixture (DPM) formulation is an extremely difficult problem. The current Monte Carlo approaches for DPM models, including the one proposed in the Appendix, are based on some sort of exploration of the space of clustering configurations of the observations. Unfortunately, the conditional expectation of the mixing distribution, given the clustering, is highly variable; much more so than the conditional expectation of the mixture density. Consequently, one needs an thorough exploration of the clustering space to obtain a reliable estimate of the mixing distribution. This is nearly impossible to achieve in finite time as this space grows exponentially with the number of observations.

4.1. *Regular mixtures.* In this subsection we will consider two regular mixture models—*regular* in the sense that $f$ is a density with respect to Lebesgue measure and smooth on its interval of support—namely, the Beta–Normal (BN) and the Gamma–Poisson (GP) mixtures

(BN)        $\theta_i \sim \frac{1}{3}\text{Beta}(3, 30) + \frac{2}{3}\text{Beta}(4, 4), \qquad X_i|\theta_i \sim N(\theta_i, \sigma^2),$

(GP)     $\theta_i \sim \text{TruncGamma}(2, 0.4)$,          $X_i | \theta_i \sim \text{Poisson}(\theta_i)$.

In each case, the samples are independent across $i = 1, \ldots, n$. Here the usual Gamma$(2, 0.4)$ distribution is truncated to $\Theta = [0, 50]$. One can easily check that conditions A2–A5 are verified for these models; in particular, A5 follows immediately from the compactness of $\Theta$. For (BN) we choose $\sigma = 0.1$ but our conclusions hold for a range of $\sigma$ containing 0.1. We also choose $f_0$ to be a Unif$(\Theta)$ density in each case.

Figures 1 and 2 display the estimates for model (BN) and (GP), respectively. In each figure, the upper left-hand cell shows the four estimates for a randomly selected run, while the other cells show the corresponding 100 estimates. The traditional estimates—MLE and NPB—of $f$ are quite poor, with the MLE being discrete and the NPB being very spiky; see Remark 7. Note that the average ESS over the 100 datasets for the (BN) and (GP) models are 538 and 324, respectively. On the other hand, RE and PARE are much more stable across samples. Moreover, as expected from the Rao–Blackwellization, we see less variability in the PARE than in the RE, in the sense that $\hat{f}_n$ hugs the true $f$ closer than does $f_n$.

That the sampling distribution is discrete in model (GP) has an interesting implication. In Figure 2 there is a (false) peak at zero for the mixture density $p_n$. This is due to the fact that the data $X_1, \ldots, X_n$ were generated by replicating each value according to its count. That is, the data sequence consists of all the 0s first, followed by all the 1s, etc. Therefore, permutation is necessary for count data stored in a particular deterministic order.

Table 1 displays the mean computation time (in seconds). In each case, the computation time for PARE is significantly less than that of NPB. In the Beta–Normal example, PARE is more than 10 times faster than the MLE, but the latter is only slightly more efficient in the Gamma–Poisson example. One explanation for this discrepancy is that PARE must process each $X_i$ individually, whereas the MLE allows for a reduction to the frequency table $n_x = \#\{i : X_i = x\}$, which can result in a significant decrease in computation time, especially when the sample size $n$ is large.

Figure 3 summarizes the $L_1$ distances $L_1(p, \hat{p})$ (left) as well as what we call a *bias-spread* summary (right) for the 100 estimates in the two regular examples. This bias-spread summary is similar to the traditional bias-variance decomposition of mean-square error: if $\hat{p}_{nt}$ is an estimate of $p$ based on the $t$th sample ($t = 1, \ldots, T$) of size $n$, then

$$
(8) \quad \text{Bias} = \int_{\mathcal{X}} |\hat{p}_{n\cdot} - p| \, d\nu \quad \text{and} \quad \text{Spread} = \frac{1}{T} \sum_{t=1}^{T} \int_{\mathcal{X}} |\hat{p}_{nt} - \hat{p}_{n\cdot}| \, d\nu,
$$

where $\hat{p}_{n\cdot}(x) = T^{-1} \sum_{t=1}^{T} \hat{p}_{nt}(x)$ is the point-wise average of the $T$ estimates of $p(x)$. We consider the sum of the bias and spread as a measure of overall variability and look at how the two components contribute to the sum. In both examples, PARE performs better in terms of overall variability, spread and, most importantly, $L_1$ loss. Compared to the other estimates, it appears that PARE does a better job of simultaneously controlling bias and spread. In the Beta–Normal example, RE also performs well. Due to the deterministic ordering issue mentioned above, RE performs quite poorly for (GP) and is not displayed. Note that these relative comparisons remain the same when the $L_1$ distance is replaced by the KL divergence.

TABLE 1
*Mean computation time (in seconds) for PARE, MLE and NPB over the $T = 100$ samples. RE (not displayed) is about 100 times faster than PARE*

| Model | PARE | MLE | NPB |
|---|---|---|---|
| BN | 0.14 | 1.11 | 43.77 |
| GP | 0.12 | 0.20 | 31.41 |

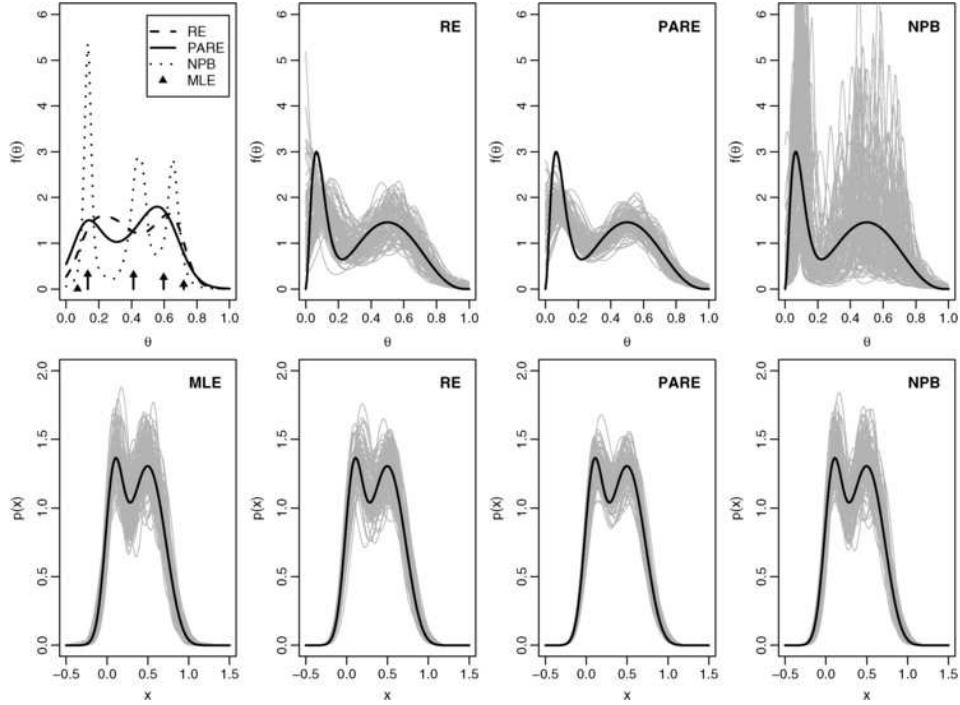

FIG. 1. *Plots of the mixing density estimates (top row) and corresponding mixture density estimates (bottom row) for model (BN). The upper left-hand cell shows all four mixing density estimates for a randomly selected run, while the remaining cells show the true f or p (black) with the $T = 100$ estimates (gray).*

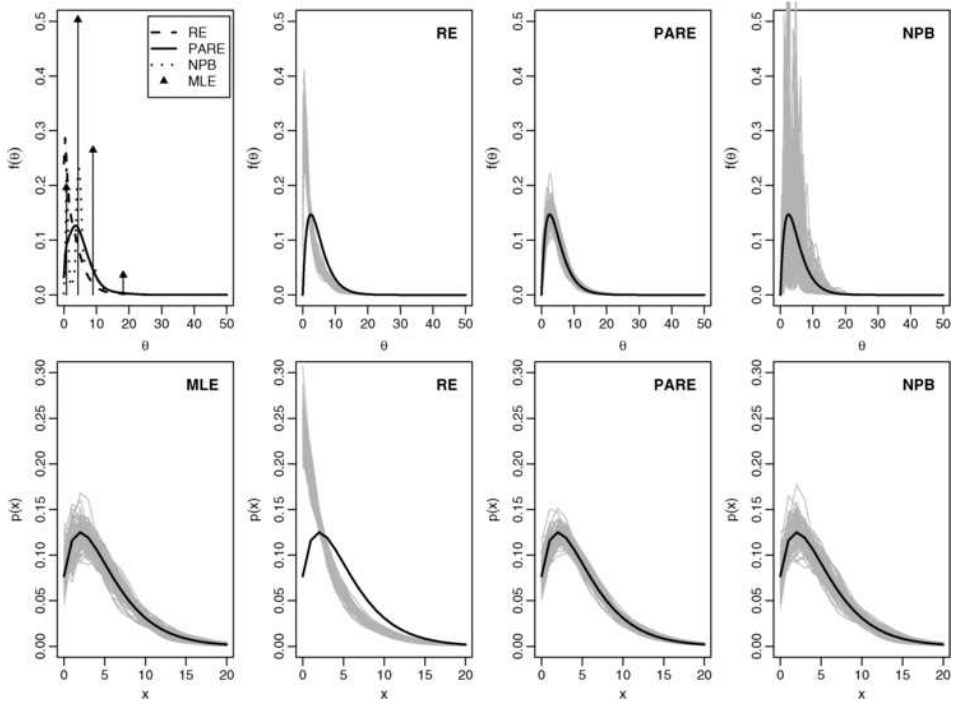

FIG. 2. *Plots of the mixing density estimates (top row) and corresponding mixture density estimates (bottom row) for model (GP). The upper left-hand cell shows all four mixing density estimates for a randomly selected run, while the remaining cells show the true $f$ or $p$ (black) with the $T = 100$ estimates (gray).*

4.2. *Irregular mixture.* For an irregular mixture, we take $f$ to have both a discrete and an absolutely continuous component. In particular, consider the Irregular–Normal (IN) mixture

(IN)    $\theta_i \sim \frac{2}{3}\delta_{\{0\}} + \frac{1}{3}\text{TruncNormal}(0, 4), \qquad X_i|\theta_i \sim N(\theta_i, 1),$

where the samples are independent across $i = 1, \ldots, n$, $\delta_{\{0\}}$ denotes a point-mass at zero, and the usual $N(0, 4)$ distribution is truncated to $\Theta = [-10, 10]$. Note that the choice of dominating measure $\mu$ is Lebesgue measure on $\Theta$ plus a unit mass at zero. The initial guess/hyperparameter $f_0$ is taken to be $\frac{1}{2}\delta_{\{0\}} + \frac{1}{2}\text{Unif}(\Theta)$ density. In this subsection we focus on just the PARE and NPB estimates, the top two performers in Section 4.1.

Figure 4 shows the 100 estimates of the the *absolutely continuous* part $f_{\text{ac}}$ of the mixing distribution as well as the corresponding estimates of the mixture. Just as in Section 4.1, we see PARE has considerably less variability than NPB (with an average ESS of about 330) on the $\theta$-scale, while both perform comparably on the $x$-scale. The left-most plot in Figure 5 summarizes the 100 estimates $\hat{\pi}$ of $\pi = \Pr(\theta = 0)$. Both procedures tend to overestimate $\pi = 0.667$ (horizontal line). Most likely, this is because $f_{ac}$ is also fairly concentrated around $\theta = 0$. The right two plots in Figure 5 summarize $L_1(p, \hat{p})$ and the bias-spread over the 100 samples. PARE, again, tends to be much more accurate under $L_1$ loss: on average, $L_1(p, p_{\mathrm{NPB}})$ is about 34% larger than $L_1(p, \hat{p}_n)$. Also, PARE seems to handle the twin bias-spread problems better than NPB.

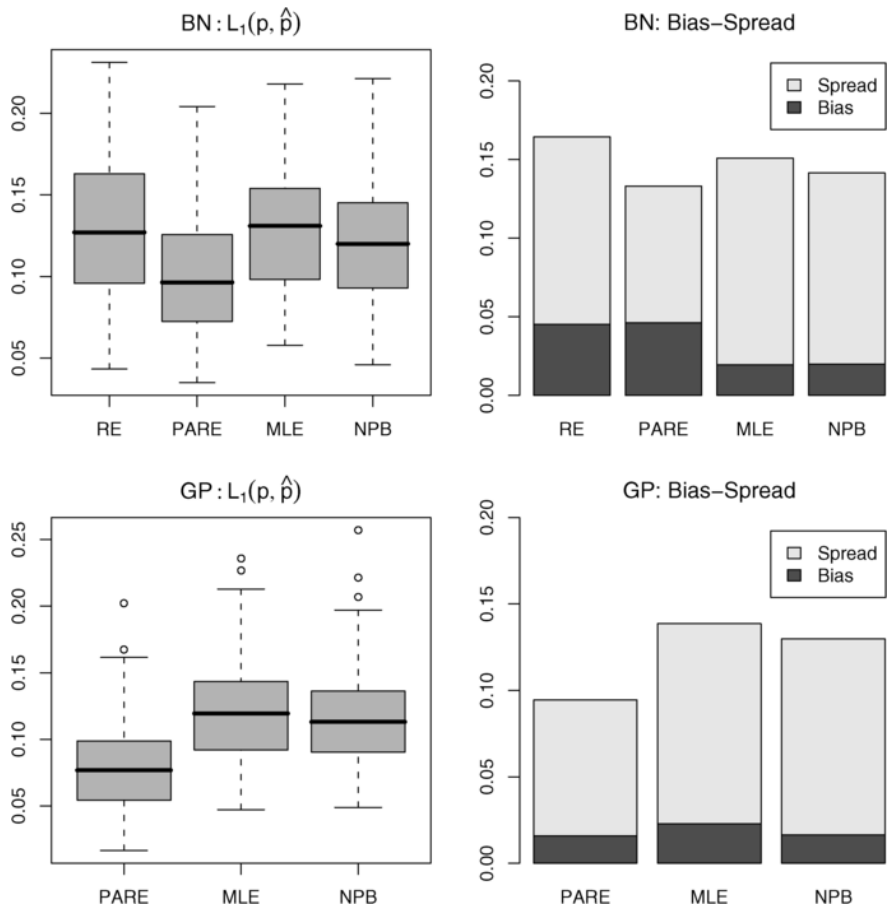

Fig. 3. *Summary of the $L_1$ distance $L_1(p, \hat{p})$ (left column) and Bias–Spread tradeoff (right column) for models (BN) (top row) and (GP) (bottom row).*

4.3. *Massive data example.* The irregular mixture (IN) in Section 4.2 arises in many important applications. In microarray analysis [29] or quantitative trait loci (QTL) mapping [3], each $\theta$ represents the expression level of a single gene or the association level of a single genetic marker, respectively. For the nonparametric regression problem [4], the $\theta$'s represent coefficients of, say, a wavelet basis expansion of the regression function. In each example, the $\theta$-vector is assumed to be *sparse* in the sense that most of the $\theta$'s are zero. To account for sparseness, a Bayesian formulation assumes that the $\theta$'s are independent observations from a common prior distribution

$$(9) \qquad F(d\theta) = \pi \delta_{\{0\}}(d\theta) + (1 - \pi) f_{ac}(\theta) \, d\theta.$$

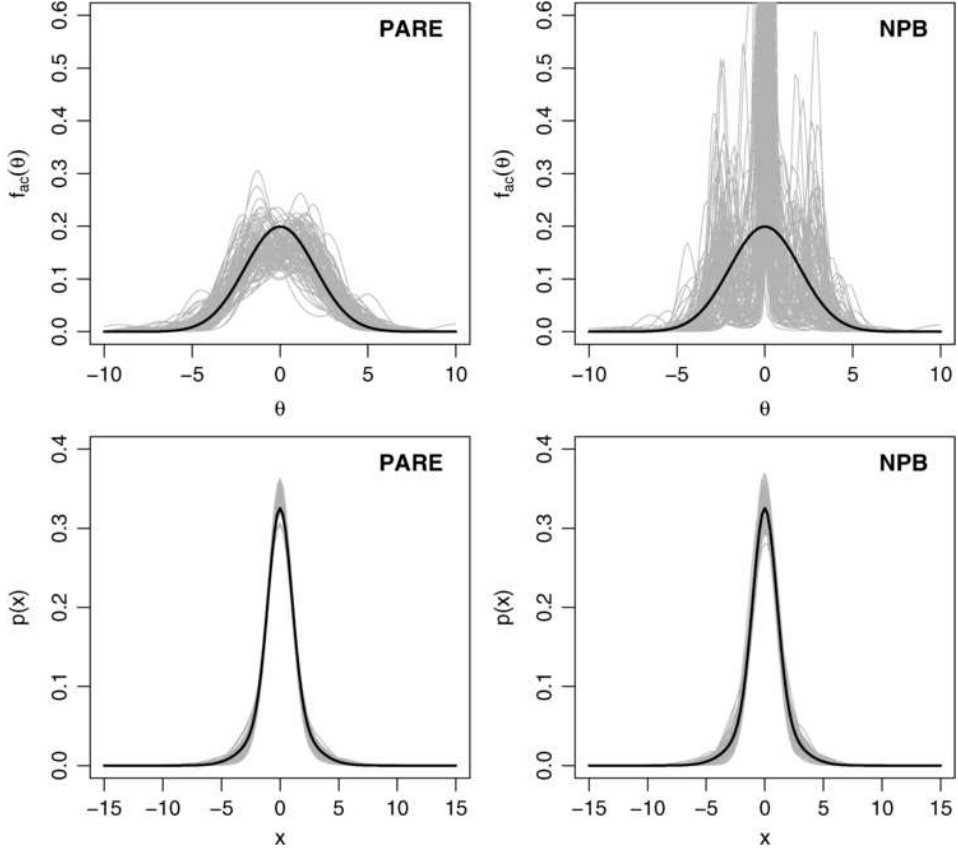

Fig. 4. *Plots of the absolutely continuous parts of the mixing distributions (top row) and corresponding mixture density estimates (bottom row) for model (IN). The true $f_{ac}$ or $p$ are shown in black with the $T = 100$ estimates in gray.*

A fully Bayesian analysis can be difficult in these applications: the results are very sensitive to the choice of hyperparameters $(\pi, f_{ac})$ [4, 29]. However, the dimension $n$ of the $\theta$-vector is quite large so an *empirical Bayes* approach [28] is a popular alternative. It was shown in Section 4.2 that both the PARE and NPB can be used to estimate $(\pi, f_{ac})$, but when $n$ is extremely large, computation becomes much more expensive, particularly for NPB.

We take a simulated dataset of size $n = 50{,}000$ from the model (IN) in Section 4.2. Figure 6 shows the PARE and NPB estimates of $(\pi, f_{ac})$ in (9).

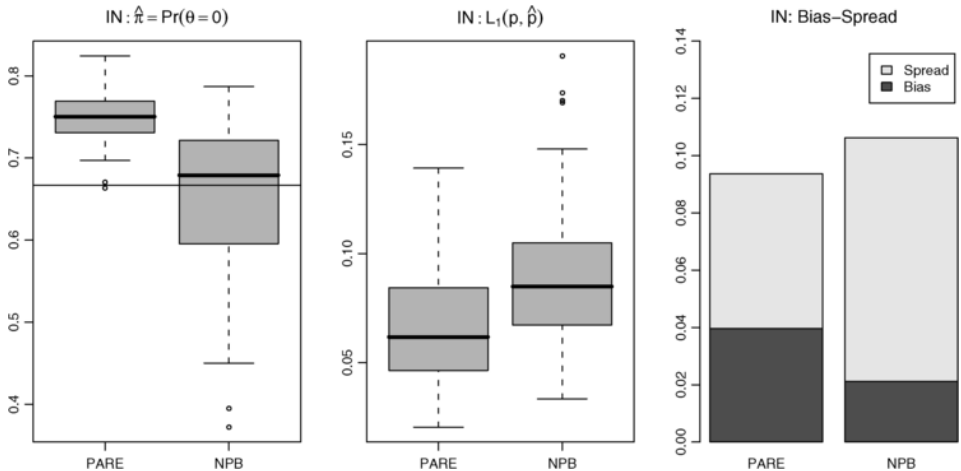

FIG. 5. *Summary of the estimates of $\pi = \mathrm{Pr}(\theta = 0)$ (left), summary of the $L_1$ distance $L_1(p, \hat{p})$ (middle) and Bias–Spread tradeoff (right) for model (IN).*

Although the PARE has two modes, it is a much closer approximation to the true $f_{\mathrm{ac}}$ compared to the spiky NPB estimate. An important point is that, even with 10,000 importance samples, *the ESS is only 1*; see Remark 7. The estimates $\hat{\pi}$ are also displayed and $\hat{\pi}_{\mathrm{PARE}} = 0.733$ and $\hat{\pi}_{\mathrm{NPB}} = 0.772$ are both slightly larger than the target $\pi = 0.667$. Figure 6 also shows the estimates $\hat{p}$ of the marginal density $p$. With $n = 50,000$, $L_1$-consistency of $p_{\mathrm{NPB}}$ [2, 11] and $\bar{p}_n$ has kicked in, and that of NPB and PARE follows by Remark 6. What is perhaps most important in massive data problems— where almost any estimate will perform well—is computational efficiency. Here, the PARE was obtained in *45 seconds*, while the NPB estimate took nearly *6 hours*.

Evidence suggests that PARE is a much better procedure than NPB in this problem for an empirical Bayes analysis. Compared to NPB, the PARE algorithm is easier to implement, the computation is significantly faster, and the resulting estimates of $(\pi, f_{\mathrm{ac}})$ are much more accurate.

**5. Discussion.** The previous analyses in [12, 20] fell short of proving strong consistency of recursive estimate $f_n$ in the general case, each only establishing convergence for the case of known and finite $\Theta$. Here a more general theorem is proved by extending the martingale approach taken by Ghosh and Tokdar [12], namely, by extending the approximate martingale representation of $K(f, f_n)$ on the $\Theta$-space to $K(p, p_n)$ on the $\mathcal{X}$-space. That the KL is the appropriate divergence measure to use for our purposes is not immediately clear, but the stochastic approximation representation of Newton's algorithm for known finite $\Theta$ along with the Lyapunov function properties of the KL divergence shown in Martin and Ghosh [20], show that the KL divergence is, indeed, quite natural. This stochastic approximation representation of the recursive algorithm in [20] continues to hold for more general $\Theta$ and we speculate that an alternative proof of convergence can be given based on this fact. Unfortunately, definitive, ready-to-use results on convergence of stochastic approximation algorithms in such general spaces are, to our knowledge, not yet available.

The failure of these previous analyses [12, 20] suggested that sample paths of the recursive estimate were, in some sense, unstable. In keeping with the stochastic approximation representation of the algorithm, we considered a *stabilized* version of $f_n$, namely,

$$f_{n:W} = \frac{\sum_{i=1}^{n} w_i f_i}{\sum_{i=1}^{n} w_i},$$

which is a weighted average of the iterates $f_i$ of the recursive algorithm. This technique of averaging the iterates, common in the stochastic approximation literature, can often improve stability properties of the algorithm, such as decreasing the variance of an estimate reached in finite time or increasing the rate of convergence; see, for example, Kushner and Yin [15]. While $f_{n:W}$ performs quite poorly compared to $f_n$ in the cases we considered, it was in proving $f_{n:W} \to f$ that Theorem 3 was discovered, opening the door to the consistency results presented in Section 2.

In simulations (including others not presented here), we have observed that $f_n$ converges quite rapidly to the true mixing density $f$. For weights of the form $w_i = (i^\alpha + 1)^{-1}$ for $\alpha \in (0.5, 1]$, the convergence was typically fastest for $\alpha = 1$. These simulation results, together with the stochastic approximation representation [20] of the recursive algorithm and the well-known results on convergence of stochastic approximation algorithms [15], suggest the following conjecture: $K_n^* = O_P(w_n^\beta)$ for some $\beta \in (0, 1)$. While the numerical evidence is consistent with this conjecture, the rate of convergence remains an open problem.

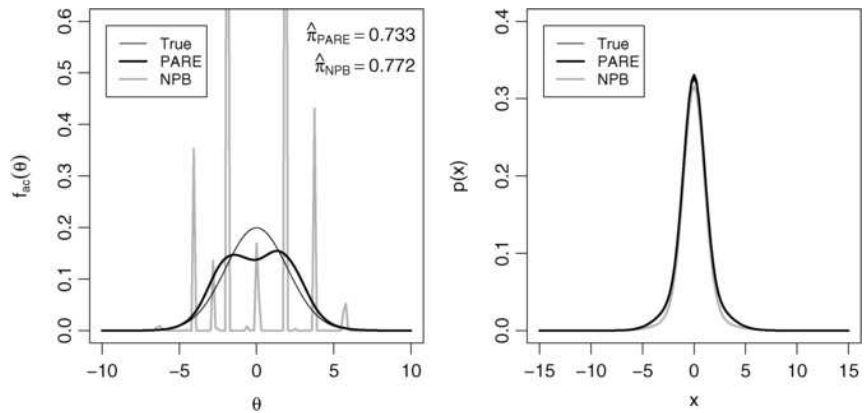

FIG. 6. *Plot of the absolutely continuous part of the mixing distribution (left) and corresponding mixture density estimates (right) for model (IN) in the massive data example.*

A drawback of the recursive algorithm is that it cannot handle an additional unknown parameter $\xi$ in $p(\cdot|\theta)$, such as an unknown $\xi = \sigma^2$ when $p(\cdot|\theta)$ is a $N(\theta, \sigma^2)$ density. Martin and Ghosh [20] tackle this problem when replicates $X_{i1}, \ldots, X_{ir}$ are available from $p(\cdot|\theta_i, \xi)$. The general idea is to use a suitable estimate $\hat{\xi}_i = \hat{\xi}(X_1, \ldots, X_i)$ of $\xi$, based on the first $i$ observations, as a plug-in in the update $f_{i-1} \mapsto f_i$ in 2. This procedure has performed well in a variety of simulations but convergence results are known only for the case of finite $\Theta$ [20]. The proof of convergence in [20] is based on a stochastic approximation representation of the algorithm and, therefore, does not easily extend to more general $\Theta$.

The numerical illustrations given here, as well as in [12, 20, 22], suggest that RE and PARE perform quite well in a variety of problems compared to other alternatives, such as MLE or NPB. While these alternatives are popular and have well-known convergence properties, which provide practical and theoretical justification for their use in applications, they lack the computational efficiency of the recursive algorithm and often produce very poor estimates. Even if one insists on a more traditional analysis, RE or PARE could be used in a computationally inexpensive preliminary analysis [26] to help choose an appropriate model to be fit to the observed data.

The theoretical results of the present paper establish the consistency properties the recursive algorithm was lacking which, combined with its generality, strong finite-sample performance and speedy computation, should put $f_n$ and $\hat{f}_n$ among the front-runners of mixing density estimates.

## APPENDIX: A NEW ALGORITHM FOR NPB

Consider the Dirichlet process mixture (DPM) model

$$x_i|\theta_i \sim p(\cdot|\theta_i), \qquad \theta_i|f \sim f, \qquad f \sim \mathscr{D}(c, f_0),$$

independently for $i = 1, \ldots, n$, where $p(x|\theta)$ is the likelihood function, the density $f$ on $\Theta$ is the parameter of interest, and $\mathscr{D}(c, f_0)$ is a Dirichlet process distribution with precision parameter $c > 0$ and base density $f_0$. Ferguson [10] shows that if the $\theta_i$'s were observed, then the posterior distribution $f$ is easily obtained. However, special techniques such as data augmentation [32], are needed when only the indirect observations $x_1, \ldots, x_n$ are available. In this section, we briefly outline a new approach to this problem.

In this approach, the mixing parameters $\theta_i$'s are collapsed onto only the clustering configuration $s = (s_1, \ldots, s_n)$, where $s_i$'s are defined sequentially as follows: $s_1 = 1$ and, for $i = 2, \ldots, n$

$$s_i = \begin{cases} s_j, & \text{if there is a } j < i \text{ such that } \theta_i = \theta_j, \\ 1 + \max_{j < i} s_j, & \text{otherwise.} \end{cases}$$

Like in Liu [19], the basic idea is to sequentially generate from

$$(10) \qquad p(s_t | x_1, s_1, \ldots, x_{t-1}, s_{t-1}, x_t), \qquad t = 1, \ldots, n$$

and calculate the importance weight

$$(11) \qquad \omega = p(x_1) \prod_{t=2}^{n} p(x_t | x_1, s_1, \ldots, x_{t-1}, s_{t-1}).$$

The current method differs from that of Liu [19] in two important ways. First, simulation of $s$ in 10 requires no advanced sampling techniques while the computational complexity of Liu's step (A) is problem-specific. Second, the conditional mean $f^{(s)} = \mathbb{E}(f | x, s)$ of the mixing density given the data $x$ and the clustering configuration $s$ can be easily calculated:

$$f^{(s)} = \frac{1}{c+n} \left[ cf_0 + \sum_{\ell=1}^{M} n_\ell f^{(\ell)} \right],$$

where $M = \max_j s_j$ is the total number of clusters, $n_\ell = \#\{j : s_j = \ell\}$ are the cluster sizes and $f^{(\ell)}$ are the cluster specific "parametric" posterior mean densities given by

$$f^{(\ell)}(\theta) \propto \prod_{j : s_j = \ell} p(x_j | \theta) f_0(\theta).$$

These calculations are summarized in the following algorithm.

1. Set $M = 1$, $s_1 = 1$, $n_1 = 1$,

$$f^{(1)}(\theta) = \frac{p(x_1 | \theta) f_0(\theta)}{\int p(x_1 | \theta') f_0(\theta') \mu(d\theta')} \quad \text{and} \quad \omega = \int p(x_1 | \theta) f_0(\theta) \mu(d\theta).$$

2. For $i = 2, \ldots, n$ repeat
   (a) Set $q_0 = c \int p(x_i | \theta) f_0(\theta) \mu(d\theta)$ and compute

   $$q_\ell = n_\ell \int p(x_i | \theta) f^{(\ell)}(\theta) \mu(d\theta), \qquad \ell = 1, \ldots, M.$$

   (b) Update $\omega \leftarrow \omega \sum_{\ell=0}^{M} q_\ell / (c + i - 1)$.
   (c) Draw $m$ from $\{0, 1, \ldots, M\}$ with probabilities $(p_0, p_1, \ldots, p_M)$ where $p_\ell \propto q_\ell$.

(d) If $m = 0$, then update $M \leftarrow M + 1$, set $s_i = M$, $n_M = 1$ and

$$f^{(M)}(\theta) = \frac{p(x_i|\theta)f_0(\theta)}{\int p(x_i|\theta')f_0(\theta')\mu(d\theta')}.$$

Otherwise, set $s_i = m$ and update $n_m \leftarrow n_m + 1$ and

$$f^{(m)}(\theta) \leftarrow \frac{p(x_i|\theta)f^{(m)}(\theta)}{\int p(x_i|\theta')f^{(m)}(\theta')\mu(d\theta')}.$$

Steps 1 and 2 are repeated $R$ times independently, producing estimates $f^{(s_r)}$ and weights $\omega_r$, for $r = 1, \ldots, R$. Then, based on the identity $f_{\mathrm{NPB}} = \mathbb{E}[f^{(s)}]$, the posterior mean is approximated by the weighted average

$$(12) \qquad f_{\mathrm{NPB}} = \frac{1}{\omega_1 + \cdots + \omega_R} \sum_{r=1}^{R} \omega_r f^{(s_r)}.$$

Note, finally, that permuting the observations $x_1, \ldots, x_n$ before each of these $R$ iterations can greatly improve the efficiency of the algorithm.

**Acknowledgments.** The authors would like to thank the Associate Editor and two referees for their insightful comments and suggestions.

S. T. TOKDAR
DEPARTMENT OF STATISTICS
CARNEGIE MELLON UNIVERSITY
BAKER HALL
PITTSBURGH, PENNSYLVANIA 15213
USA
E-MAIL: stokdar@stat.cmu.edu

R. MARTIN
DEPARTMENT OF STATISTICS
PURDUE UNIVERSITY
250 NORTH UNIVERSITY STREET
WEST LAFAYETTE, INDIANA 47907
USA
E-MAIL: martinrg@stat.purdue.edu

J. K. GHOSH
DIVISION OF THEORETICAL STATISTICS
   AND MATHEMATICS
INDIAN STATISTICAL INSTITUTE
203 B T ROAD
KOLKATA, INDIA 700108
AND
DEPARTMENT OF STATISTICS
PURDUE UNIVERSITY
250 NORTH UNIVERSITY STREET
WEST LAFAYETTE, INDIANA 47907
USA
E-MAIL: ghosh@stat.purdue.edu